\documentclass[10pt,leqno]{article}
\usepackage[OT2,OT1]{fontenc}
\usepackage[russian,english]{babel}
\usepackage{amssymb,amsmath,amsthm}
\usepackage{array}
\usepackage{float}
\usepackage{a4}

\makeatletter

\renewcommand\section{\@startsection {section}{1}{\z@}%
                                   {-3.5ex \@plus -1ex \@minus -.2ex}%
                                   {2.3ex \@plus.2ex}%
                                   {\normalfont\large\bfseries}}
\renewcommand\subsection{\@startsection{subsection}{2}{\z@}%
                                     {-3.25ex\@plus -1ex \@minus -.2ex}%
                                     {1.5ex \@plus .2ex}%
                                     {\normalfont\normalsize\bfseries}}

\newtheoremstyle{styleDef}{10pt}{10pt}{\upshape}{}{\scshape}{.}{7pt}{\thmnumber{{\upshape\bfseries#2.}\ }\thmname{#1}\thmnote{ #3}}
\newtheoremstyle{styleThm}{10pt}{10pt}{\slshape}{}{\scshape}{.}{7pt}{\thmnumber{{\upshape\bfseries#2.}\ }\thmname{#1}\thmnote{ #3}}
\newtheoremstyle{styleRem}{10pt}{10pt}{\upshape}{}{\itshape}{.}{7pt}{\thmnumber{{\upshape\bfseries#2.}\ }\thmname{#1}\thmnote{ #3}}

\theoremstyle{styleDef}
  \newtheorem{Definition}{Definition}

\theoremstyle{styleThm}
  \newtheorem{Theorem}[Definition]{Theorem}
  \newtheorem*{Theorem*}{Theorem}
  \newtheorem*{Statement*}{Statement}
  \newtheorem*{Question*}{Question}
  \newtheorem{Proposition}[Definition]{Proposition}
  \newtheorem*{Proposition*}{Proposition}
  
  \newtheorem{Lemma}[Definition]{Lemma}

\theoremstyle{styleRem}
  \newtheorem*{Remark}{Remark}

\newenvironment{myitemize}{\list{$\bullet$}{\itemsep=0pt}}{\endlist}  

\newcommand*\tabstrut{\vrule depth 6pt height 13pt width 0pt}

\newcommand*\sseq{}

\let\sseq\subseteq                  
\let\le\leqslant \let\leq\leqslant  
\let\ge\geqslant \let\geq\geqslant  

\newcommand*\scal[2]{\langle #1\mid\nobreak #2\rangle} 

\newcommand*\card[1]{\left|#1\right|}      
\newcommand*\Har[1]{\mathcal{H}^{(#1)}}    
\newcommand*\Pol[1]{\mathcal{P}^{(#1)}}    

\newcommand*\NN{\mathbb{N}}   
\newcommand*\ZZ{\mathbb{Z}}   
\newcommand*\RR{\mathbb{R}}   
\newcommand*\CC{\mathbb{C}}   

\newcommand*\Se{\mathbb{S}}   

\newcommand*\Orth{\mathrm{O}} 
\newcommand*\SL{\mathrm{SL}}  
\newcommand*\Aut{\mathop{\mathrm{Aut}}} 

\newcommand*\HH{\mathbb{H}}   

\newcommand*\bound{\mathrm{b}}
\newcommand*\boundLP{\bound_\mathrm{LP}}
\newcommand*\boundD{\bound_\mathrm{D}}
\newcommand*\boundY{\bound_\mathrm{Y}}
\newcommand*\Bound{\mathrm{B}}

\newcommand*\Qp[1]{\mathrm{Q}^{(#1)}}

\newcommand*\smatr{\left(\smallmatrix}
\newcommand*\esmatr{\endsmallmatrix\right)}
\newcommand*\smatrabcd{\smatr a&b\\c&d\esmatr}

\newcommand*\legendre{\genfrac(){}{}}  

\newcommand*\Weyl{\mathrm{W}}          

\newcommand*\A{\mathbf{A}}
\newcommand*\D{\mathbf{D}}
\newcommand*\E{\mathbf{E}}
\newcommand*\F{\mathbf{F}}

\newcommand*\BW{\mathit{BW}} 
\newcommand*\KZ{\varGamma}   

\newcommand*\Mod{\mathcal{M}}               
\newcommand*\Pod{{}^{0}\!\mathcal{M}}       
\newcommand*\varMod{\mathcal{\widehat M}}         
\newcommand*\varPod{{}^{0}\!\mathcal{\widehat M}} 
\newcommand*\Lat{\mathcal{L}\!\mathit{at}}  

\DeclareMathOperator\sym{sym}  
\newcommand*\Sym{\mathfrak{S}} 

\title{Construction of spherical cubature formulas using lattices}
\author{Pierre \textsc{de la Harpe}, Claude \textsc{Pache}, Boris \textsc{Venkov}%
        \footnote{The authors acknowledge support from the \emph{Swiss National Science Foundation}.\endgraf
                  Pierre de la Harpe and Claude Pache,
                  Section de Math\'ematiques, Universit\'e de Gen\`eve, C.P.~64, 1211~Gen\`eve~4, Switzerland.\endgraf
                  Boris Venkov, Petersburg Department of Steklov Institute of Mathematics, Fontanka~27,
                  191011~St.~Petersburg, Russia.\endgraf
                  Email: Pierre.delaHarpe@math.unige.ch, Claude.Pache@math.unige.ch, bbvenkov@yahoo.com
                 }\\
        {\small Contact author: Claude Pache (Claude.Pache@math.unige.ch)}%
       }
\date{3 November 2005}

\begin{document}
\maketitle

\begin{abstract}
  We construct cubature formulas on spheres
  supported by homothetic images of shells in some Euclidian lattices.
  Our analysis of these cubature formulas uses results
  from the theory of modular forms.
  Examples are worked out on $\Se^{n-1}$
  for $n=4$,~$8$, $12$, $14$, $16$, $20$, $23$, and~$24$,
  and the sizes of the cubature formulas we obtain
  are compared with the lower bounds given by Linear Programming.
\end{abstract}


\section{Introduction}\label{SectIntro}

   For a dimension $n \ge 2$ and a positive number $r$, 
let $\Se^{n-1}_r$ denote the sphere of equation
$\scal{x}{x} = r$ in the Euclidean space $\RR^n$ 
given with the canonical scalar product 
$\scal{\cdot }{\cdot }$,
and let $\sigma$ denote the rotation-invariant measure on such a sphere,
normalised by $\sigma(\Se^{n-1}_r) = 1$; 
we write $\Se^{n-1}$ for $\Se^{n-1}_1$.
For an integer $t \ge 0$,
a \emph{cubature formula of strength $t$ on $\Se^{n-1}_r$}
is a pair $(X,W)$, where $X$ is a finite subset of $\Se^{n-1}_r$
and where $W : \Se^{n-1}_r \to \RR_{>0}$
is a positive-valued function such that
\[
\sum_{x \in X} W(x)\,f(x)  = 
\int_{\Se^{n-1}_r} f(x)\,d\sigma(x)
\tag{CF}
\]
for every polynomial function $f:\RR^n\to\CC$
of degree at most~$t$.
The set $X$ is the \emph{support} of the cubature formula $(X,W)$,
the points in $X$ are its \emph{nodes},
the values $W(x)$ are its \emph{weights}, 
and the cardinality $\card{X}$ of $X$ is its \emph{size}.
A \emph{spherical $t$-design} on $\Se^{n-1}_r$
is a finite subset $X$ of $\Se^{n-1}_r$ which,
together with the constant weight
(by necessity of value $\card{X}^{-1}$),
is a cubature formula of strength~$t$.

   There are existence results according to which,
for any dimension $n$ and strength $t$,
there exist cubature formulas, and indeed spherical designs,
of strength $t$ on $\Se^{n-1}$.
For cubature formulas, there are elementary arguments
using general methods of convexity,
which give some upper bounds on sizes
(see Section~\ref{SectBounds});
the existence of spherical designs is
a particular case of results of \cite{SeymourZaslavsky},
which do not provide meaningful size bounds.
There are no known straightforward constructions,
except for $n=2$, 
in which case, for any $t \ge 0$,
the $(t+1)$th roots of unity in $\CC$
provide a spherical $t$-design on the circle;
and in larger dimensions for low values of $t$
(antipodal pairs, regular simplices, and regular hyperoctahedra
are respectively spherical $1$-, $2$-, and $3$-designs,
of the lowest possible sizes).

  It is therefore a natural question 
to ask for more explicit constructions.
Some have been given in terms of orbits of finite groups
on the sphere (see \cite{HarPac}),
and by other various methods
(see \cite{Bajnok} and later papers by Bajnok,
\cite{Kuperberg1}, \cite{Kuperberg2}, and \cite{HardinSloane96} for $\Se^2$).

   The object of the present paper is to describe constructions involving 
\emph{shells (or layers) of integral lattices} in $\RR^n$.
They are elaborations of constructions which can also be found in
\cite{MartinetVenkov}, \cite{hpv}, and \cite{Gaeta}.

   Let $E$ be a discrete subset of $\RR^n$.
For any $m > 0$, denote by $E_m$ the \emph{shell} $E \cap \Se^{n-1}_m$.
If $m_1, \hdots, m_r$ are pairwise distinct 
positive numbers
such that each shell $E_{m_j}$ is non-empty,
the union of the sets $\frac{1}{\sqrt{m_j}}E_{m_j}$ 
is a finite subset (possibly with multiplicities) 
of the unit sphere in $\RR^n$.
Spherical designs and cubature formulas which appear below
have supports of this kind,
where $E$ is either an integral lattice in $\RR^n$
or a union of lattices;
observe that they are all \emph{antipodal},
namely that $X=-X$ and $W(x)=W(-x)$ for all $x\in X$.

\paragraph{\boldmath Results on $\Se^3$}

   Let us describe some examples provided by this method 
when $n=4$, and compare them to known lower bounds for sizes
and to other known examples.
For $n \ge 4$,
the \emph{root lattice of type $\D_n$} is the integral lattice
\[
  D_n  = \{ x \in \ZZ^n \mid x_1 + \cdots + x_n \equiv 0 \pmod{2} \} .
\]
(We denote \emph{root systems} as $\D_n$ by bold letters
and \emph{lattices} as $D_n$ by thin letters.)
It is important for our computations that $D_4$
is a $2$-modular lattice
(the definition is recalled in Section~\ref{SectModular} below).

For $n \ge 5$, $n \ne 7$, 
it is known that each shell $(D_n)_{2m}$
is a spherical $3$-design 
which is in general not of strength $4$;
for $n = 4$, each shell $(D_4)_{2m}$
is a spherical $5$-design
which is in general not of strength $6$
(it has been checked in \cite{hpv}
that it is never of strength $6$ for $2m \le 1200$).
The sizes of the shells of the root lattice of type $\D_4$ are given
by the coefficients of the theta series
\begin{align*}
\Theta_{D_4}  
     &=  \sum_{m\ge0} \card{(D_4)_{2m}} q^{2m}
      =  \sum_{m\ge0} r_4(2m) q^{2m} \\
     &=  1 + 24\,q^2 + 24\,q^4 + 96\,q^6 + 24\,q^8 + 144\,q^{10} + 
             96\,q^{12} + 192\,q^{14} + \cdots
\end{align*}
where $r_4(2m)$ denotes the number of ways 
of writing $2m$ as a sum of four squares.

The dual $D_4^* = \{ y \in \RR^n \mid \scal{y}{D_4} \subset \ZZ \} $
is similar to $D_4$ and the renormalized dual $D_4' = \sqrt{2}\,D_4^*$
is isometric to $D_4$; 
moreover the intersection $D_4 \cap D_4'$ is reduced to $\{0\}$,
since the non-zero coordinates of the vertices in $D_4'$
are never in $\ZZ$.
For our analysis, we use some invariant theory
of the finite group $\Aut(D_4\cup D_4')$.

\begin{Theorem*}[$D_4$] Set $E = D_4 \cup D_4'$
\begin{enumerate}
\item For any even integer $2m \ge 2$, the shell $E_{2m}$ is a 
  spherical $7$-design; in particular, $E_2$ is a spherical $7$-design
  of size $48$ in $\Se^3_2$.
\item The two shells $E_2$, $E_6$ provide the support
  $X = \frac{1}{\sqrt 2}E_2 \cup \frac{1}{\sqrt 6}E_6$
  of a cubature formula of strength $11$ and size $240$.
\item Similarly, 
  $\frac{1}{\sqrt 2}E_2 \cup \frac{1}{6}E_6 \cup \frac{1}{\sqrt{10}}E_{10}$
  is the support of a cubature formula of strength~$15$ and size $528$.
\end{enumerate}
\end{Theorem*}

Claim (i) for $E_2$ appears already in Section 6 of \cite{GoethalsSeidel81},
where the authors  observe that $E_2$ 
is an orbit of the Weyl group of type $\F_4$. 
Claim (i) of Theorem $E_8$ below
and the analogous result involving the two shortest shells of a Leech
lattice appear also in \cite{GoethalsSeidel81}.

\medskip

For comparison, here are some known results concerning the $3$-sphere.

Dickson reports on a formula due to Liouville (1859)
which essentially shows that the system of roots $(D_4)_2$ of type $\D_4$
is a spherical $5$-design,
a formula due to Kempner (1912)
which essentially show that $(D_4\cup D_4')_2$
is a spherical design
of strength $7$,
and a formula due to Schur (1909)
related to a cubature formula
of strength $11$ and size $144$.
See \cite[pp.~717--724]{Dickson}.

There is a formula of strengh $7$ and size $46$ 
described in \cite{HardinSloane}.

There is an exceptional regular polytope in dimension four
known as the $600$-cell.
It has $120$ vertices which constitute a spherical $11$-design;
it is an orbit of a Coxeter group of type~$\mathbf{H}_4$.
Any orbit of this group is a spherical design of strength $11$,
but there is a special orbit of size $1440$ which is of strength $19$;
see \cite[Section~5]{GoethalsSeidel81}.
Another cubature formula of the same strength, $19$,
and of size $720$ appears in \cite{Salihov}:
its support is the union of the $120$ vertices and of the $600$ centers of faces
of the $600$-cell.

The best lower bounds known to us for the size of 
cubature formulas in $\Se^3$ are listed
in Tables \ref{tabCF4}~and~\ref{tabCF4b} of Section~\ref{SectResults}.

\paragraph{\boldmath Results on $\Se^7$}

   Let $\KZ_8$ be a Korkine-Zolotareff lattice,
namely an even unimodular lattice of dimension~$8$
(uniquely defined up to isometry by these properties).
Such a lattice is generated by a root system of type $\E_8$.
The sizes of its shells are given by the coefficients of the theta series
\begin{align*}
  \Theta_{\KZ_8}  
      &=  \sum_{m\ge0} \card{(\KZ_8)_{2m}} q^{2m}
       =  240 \sum_{m\ge0} \sigma_3(m) q^{2m} \\
      &=  1 + 240\,q^2 + 2160\,q^4 + 6720\,q^6 + 17520\,q^8 +  \cdots
\end{align*}
where $\sigma_3(m) = \sum_{d \vert m}d^3$.

   It is known that each shell $(\KZ_8)_{2m}$ is a spherical $7$-design
which is in general not of strength $8$.
Indeed, the shell $(\KZ_8)_{2m}$ would be a design of strength $8$
if and only if the Ramanujan function $\tau$ would vanish at $m$;
it is a conjecture of Lehmer \cite{Lehmer} 
that $\tau(m) \ne 0$ for all $m \ge 1$, 
and this has been checked for $m \le 10^{15}$ \cite{Serre};
more on this in \cite{hpv} and in Item 1.16 of \cite{Gaeta}.

\begin{Theorem*}[$E_8$]
  Set $E = \KZ_8$.
  \begin{enumerate}
    \item The two shells $E_2$, $E_4$ provide the support
      $X = \frac{1}{\sqrt 2}E_2 \cup \frac{1}{2}E_4$
      of a cubature formula of strength $11$ and size $2400$.
    \item Similarly, 
      $\frac{1}{\sqrt 2}E_2 \cup \frac{1}{\sqrt 6}E_6 
      \cup \frac{1}{2\sqrt 2}\,E_8$ 
      is the support of a cubature formula of strength $13$ and size $24 240$.
    \item Similarly, 
      $\frac{1}{\sqrt 2}E_2 \cup \frac{1}{2}E_4 
      \cup \frac{1}{\sqrt 6}E_6  \cup \frac{1}{2\sqrt 2}\,E_8$ 
      is the support of a cubature formula of strength $15$ and size $26 400$.
  \end{enumerate}
\end{Theorem*}

   The largest known lower bounds for the sizes of cubature formulas on
$\Se^7$ are given in Table~\ref{tabCF8} of Section~\ref{SectResults}.

\paragraph{\boldmath Dimensions $12$, $14$, $16$, $20$, $23$, and $24$}

   There are similar constructions on $\Se^{23}$ 
in terms of a Leech lattice, 
an even unimodular lattice $\Lambda$ of dimension $24$ 
with empty shell $\Lambda_2$ 
(uniquely defined up to isometry be these properties).
   We also describe constructions with shells of
\begin{myitemize}
  \item the $3$-modular Coxeter-Todd lattice $K_{12}$, 
  \item the $3$-modular Quebbemann lattice $Q_{14}$,
  \item the $2$-modular Barnes-Wall lattice $\BW_{16}$, 
  \item the $2$-modular Nebe lattices $N_{20}$ (three of them), 
  \item the $3$-modular Nebe lattice $N_{24}$, and
  \item the shorter Leech lattice $O_{23}$
    (which is up to isometry the unique unimodular integral lattice $\Lambda$
    in $\RR^{23}$ with $\Lambda_1 = \Lambda_2 = \emptyset$).
\end{myitemize}
Constructions with one shell of these lattices are already
described in \cite{BachocVenkov}.
Our results appear in Tables~\ref{tabCF12} to \ref{tabCF24} of Section~\ref{SectResults}.

\section{Bounds on the size of cubature formulas}\label{SectBounds}

For fixed values of~$n$ and~$t$, we are interested in constructing
cubature formulas on $\Se^{n-1}_r$ of strength~$t$
of small size; it is therefore useful to estimate the minimal size
that such a cubature formula can have.
We recall in this paragraph some general results on this
subject.

Let $\bound(n,t)$ be the minimal size of
a spherical cubature formula on $\Se^{n-1}_r$ of strength~$t$.

\paragraph{Linear Programming bound \cite[Theorem~5.10]{DelsarteGoethalsSeidel}, see also \cite{Yudin}}
Let $\Qp{k}$, $k\in\ZZ_{\ge0}$, be the orthogonal polynomials
for the scalar product
\[
  \scal{f}{g} := \int_{-1}^{1} f(t)\,g(t)\,(1-t^2)^{(n-3)/2}\,dt
\]
normalized to $\scal{\Qp{k}}{\Qp{k}} = \Qp{k}(1)$.
For a continuous function $F:[-1,1]\to\CC$, we define the numbers $F_k\in\CC$ by
\[
  F_k = \frac{\scal{F}{\Qp{k}}}{\Qp{k}(1)},
\]
so that $F = \sum_{k\ge0} F_k\,\Qp{k}$.
Set
\[
  M(n,t):= \biggl\{ F:[-1,1]\to\RR \text{ continuous}\biggm|
         \vcenter{\hbox{\strut$F(u)\ge0$ for $u\in[-1,1]$,}%
                  \hbox{\strut$F\ne0$, and $F_k\le0$ for $k>t$}%
                 }
            \biggr\};
\]
then $\bound(n,t)\ge F(1)/F_0$ for every $F\in M(n,t)$.
In other words
\[
  \boundLP(n,t):=\sup_{F\in M(n,t)} \frac{F(1)}{F_0}
\]
is a lower bound for $\bound(n,t)$.
Moreover, it is a consequence
of \cite[Therorems~I and~2.7]{NikovaNikov}
that there exists a \emph{polynomial} function $F\in M(n,t)$
for which $\boundLP(n,t) = F(1)/F_0$.

When $t$ is odd, we show:

\begin{Proposition}\label{SymmetricLinProg}
  Let $t=2s+1$, and let 
  \[
    \widehat{M}(n,t):= \{ F\in M(n,t) \mid F(-u)=F(u) \}.
  \]
  Then
  \[
    \boundLP(n,t)=\sup_{F\in \widehat{M}(n,t)} \frac{2\,F(1)}{F_0}.
  \]
\end{Proposition}

\begin{proof}
  Let $F\in M(n,t)$; set $\widehat{F}(u)=\bigl(F(u)+F(-u)\bigr)\big/2$;
  we have $\widehat{F}\in\widehat M(n,t)$.
  Set $F^*(u)=(1+u)\,\widehat{F}(u)$. One checks that
  \[
    F^*\in M(n,t).
  \]
  Yet $2\,\widehat{F}(1) = F(1)+F(-1)\ge F(1)$ and $F_0 = \widehat{F}_0 = F^*_0$,
  therefore
  \[
     \frac{F(1)}{F_0} \le \frac{2\, \widehat{F}(1)}{\widehat{F}_0} = \frac{F^*(1)}{F^*_0}.
  \]
  The result follows.
\end{proof}
  
\smallskip
In general, the exact value of $\boundLP(n,t)$ is not known.
The following bounds are obtained by chosing a particular $F$ in $M(n,t)$:

\paragraph{Delsarte (or Fisher-type) bound \cite[Theorems 5.11~and~5.12]{DelsarteGoethalsSeidel}}
Choose $F(u):=\bigl(\sum_{k=0}^{s}\Qp{k}(u)\bigr)^2$ if $t=2s$,
and $F(u):=\bigl(1+u\bigr)\bigl(\sum_{j=0}^{[s/2]}\Qp{s-2j}(u)\bigr)^2$
if $t=2s+1$. We obtain the bounds
\begin{gather*}
  \boundD(n,2s)   = \binom{n+s-1}{n-1}+\binom{n+s-2}{n-1}, \\
  \boundD(n,2s+1) = 2\,\binom{n+s-1}{n-1}.
\end{gather*}
Moreover, it is known
that, if $n\ge3$, the equality $\bound(n,t)=\boundD(n,t)$
is possible only for some values of $(n,t)$: see \cite{BannaiDamerell1}, \cite{BannaiDamerell2},
and \cite{BannaiMunemasaVenkov}.
There also exists a criterion to decide
when $\boundD(n,t) = \boundLP(n,t)$: see \cite[Theorem~I]{NikovaNikov}.

\paragraph{Yudin bound \cite{Yudin}}
Another choice of~$F$ gives the bound
\[
  \boundY(n,t)   = \frac{ \int_{-1}^1 (1-u^2)^{(n-3)/2} du }%
                        { \int_\gamma^1 (1-u^2)^{(n-3)/2} du },
\]
where $\gamma$ is the largest root of the polynomial $\bigl(\Qp{t+1}\bigr)'$.

\paragraph{Special cases} For some values of $(n,t)$,
it is possible to compute the exact value of $\bound(n,t)$.
For example, it is shown in
\cite{Andreev} that $\bound(4,11)=\boundLP(4,11)=120$.

\paragraph{Numerical estimate of the Linear Programming bound}
For fixed values of $(n,t)$, $t=2s+1$, the following procedure gives an estimate
of the Linear Programming bound $\boundLP(n,t)$.
(There is a similar procedure when $t$ is even.)
\begin{itemize}
  \item Choose a degree~$d \ge s$,
    and a finite subset $A\sseq[0,1]$ of well-distributed points,
    for example $A=\{i/N\mid i=0, 1, \dots, N\}$ with $N$ large.
  \item By linear programming, find,
    among all polynomial~$F$ satisfying
    \[
      F = 1 + \sum_{i=1}^{d} F_k \Qp{2k}, \quad
      F(u)\ge0 \text{ for $u\in A$}, \quad
      F_k\le0 \text{ for $k>t$},
    \]
    the polynomial minimizing $F(1)$;
    denote it by $G$.
  \item Set $\epsilon:= - \inf_{x\in[0,1]} G(x) \ge 0$. The polynomial
    $\widetilde{G}:= G + \epsilon$ is in the set $\widehat{M}(n,t)$,
    and we have the estimate
    \[      
      \frac{2\,\widetilde{G}(1)}{\widetilde{G}_0} = \frac{2\,\bigl(G(1)+\epsilon\bigr)}{1+\epsilon}
      \le \boundLP(n,t).
    \]
    Note that we have also $\boundLP^{(d)}(n,t)\le 2\,G(1)$,
    where
    \[
      \boundLP^{(d)}(n,t):=\sup_{\substack{F\in M(n,t)\\\deg f\le d}} \frac{F(1)}{F_0}.
    \]    
\end{itemize}
Yet, if $d$ is sufficiently large and $A$ sufficiently dense,
$2\,\widetilde{G}(1)/\widetilde{G}_0$ is a good approximation of $\boundLP(n,t)$.
In practice, we first apply the procedure with a relatively large $d$ and a relatively small set~$A$;
we observe that $G_k=0$ for $k>d_0$. Then we apply the procedure with $d=d_0$
and with a larger set~$A$.

The procedure just described imitates the one used in \cite[Chap.~13]{ConwaySloane}
to compute bounds for kissing numbers.
\bigskip

We end this section with another kind of bound.
Set
\[
  \Bound(n,t):= \binom{n+t-1}{n-1}+\binom{n+t-2}{n-1}.
\]
Note that $\Bound(n,t)=\boundD(n,2t)$. This number
is the dimension of the space of
the restrictions to the sphere of the polynomial functions of degree at most~$t$.

\begin{Proposition}
  For every $n$ and $t$, there exist cubature formulas with
  at most $\max\{1,\Bound(n,t)-1\}$ nodes.
\end{Proposition}

\begin{Proposition}
  Let $(X,W)$ be a cubature formula on~$\Se^{n-1}_r$ of strength~$t$.
  Then there exist a subset $X'\sseq X$ and a weight function $W':X'\to\RR_{>0}$
  such that $\card{X'}\le \Bound(n,t)$ and $(X',W')$ is a cubature formula of degree~$t$.
\end{Proposition}

See \cite[Proposition~2.6]{Gaeta}. The second proposition is a consequence of the proof of the first one
and provides an algorithm for computing $(X',W')$ from $(X,W)$.

Note however that these propositions do not hold for spherical designs
(case of constant weights).


\section{Harmonic polynomials}\label{SectHarmonicGroup}

It is sufficient to check condition (CF) 
in the definition of cubature formulas 
for harmonic polynomials only. 
Recall that a smooth function $f$ on $\RR^n$ is \emph{harmonic}
if $\Delta f = 0$, 
where $\Delta = \sum_{i=1}^n (\partial / \partial x_i)^2$ 
is the usual Laplacian.
We denote by $\Har{k}(\RR^n)$
the space of harmonic polynomial functions 
$\RR^n \to \CC$ that are homogenous of degree $k$.
It is a classical fact that 
the space of restrictions to $\Se^{n-1}$ of
homogenous polynomials of degree at most $t$
coincides with the space of restrictions of 
the direct sum $\bigoplus_{k=0}^t \Har{k}(\RR^n)$.
This implies the following well-known criterion;
see for example \cite[Theorem~3.2]{VenkovMartinet1}
(details are written there for spherical designs,
but the proof carries over to cubature formulas). 

\begin{Proposition}\label{CharactHarmonic}
Let $X$ be a finite subset of some sphere $\Se^{n-1}_m$ in $\RR^n$,
let $W:X \to \RR_{>0}$ be a weight function,
and let $t \ge 0$ be an integer.
Then $(X,W)$ is a cubature formula of strenght $t$ on $\Se^{n-1}_m$
if and only if $\sum_{x\in X}W(x)=1$ and
\[
  \sum_{x \in X} W(x)\,P(x)  =  0
  \qquad \text{for all $P \in \Har{k}(\RR^n)$, $1 \le k \le t$.}
\]
\end{Proposition}

Let $\Orth(n)$ be the orthogonal group of the Euclidean space $\RR^n$;
this group acts by isometries on the sphere $\Se^{n-1}_r$.
Moreover, $\Orth(n)$ acts naturally on $\Har{k}(\RR^n)$
by $g\cdot f=f\circ g^{-1}$.

Whenever a group $G$ acts on a space $V$, we denote by $V^G$
the subspace of all elements of~$V$ fixed by~$G$.

Proposition~\ref{CharactHarmonic} above can be refined in the case
the cubature formula is invariant by a finite subgroup of~$G$.
Namely, we have the following result, known as \emph{Sobolev's theorem}
(\cite{Sobolev}, see also \cite[Chap.~2, \S~2, Theorem~2.3]{SobolevBook}):

\begin{Proposition}\label{CharactHarmonicInv}
  Let $X$ be a finite subset of $\Se^{n-1}_r$ and $W:X\to\RR_{>0}$
  be a weight function. Let $G$ be a finite subgroup of $\Orth(n)$ leaving
  $(X,W)$ invariant, that is
  \[
     gx\in X \text{ and }  W(gx)=W(x)\qquad \forall x\in X, \forall g\in G.
  \]
  Then $(X,W)$ is a cubature formula of strength~$t$ if and only if
  \[
    \sum_{x\in X} W(x)\,P(x) = 0 \qquad \forall P\in\Har{k}(\RR^n)^G, 1\le k\le t.
  \]
\end{Proposition}

\begin{proof}
  For $P\in\Har{k}(\RR^n)$, let $P^G:= \card{G}^{-1}\sum_{g\in G} g\cdot f$.
  On the one hand, the $G$-invariance of $(X,W)$ implies
  \[
    \sum_{x\in X} W(x)\,P^G(x) = \sum_{x\in X} W(x)\,P(x),
  \]
  and, on the other hand, the $G$-invariance of the Lebesgue measure implies
  \[
    \int_{\Se_r^{n-1}}P^G(x)\,d\sigma(x) = \int_{\Se_r^{n-1}} P(x)\,d\sigma(x).
  \]
  The result follows from Proposition~\ref{CharactHarmonic}.
\end{proof}

\section{Modular forms}\label{SectModular}

This section is a reminder from \cite{BachocVenkov}.
Consider a lattice $\Lambda$ in $\RR^n$
and a homogeneous harmonic polyomial $P \in \Har{k}(\RR^n)$ of degree $k$.
The \emph{theta series of $\Lambda$  with harmonic coefficient $P$}
is the formal power series $\Theta_{\Lambda,P}$ defined by
\[
  \Theta_{\Lambda,P} := \sum_{x\in\Lambda} P(x)\,q^{\scal xx}
  = P(0) + \sum_{m>0}\biggl(\sum_{x\in\Lambda_m} P(x)\biggr) q^m.
\]
In case $P=1$, we write simply $\Theta_{\Lambda}$.
Set
\[
  q := e^{i\pi z},\quad z\in\HH,
\]
with $\HH=\{z\in \CC \mid \Im z>0 \}$ the Poincar\'e half-plane;
the theta series $\Theta_{\Lambda,P}$ converges uniformly
on every compact subset of $\HH$ and defines consequently
a holomorphic function on $\HH$.
When the lattice $\Lambda$ satisfies appropriate conditions,
this function is a \emph{modular form} of \emph{weight~$\omega$},
and we can use results
on modular forms for the computation of the theta series.

Let us first describe the class of lattices which plays the most important role below.
A lattice $\Lambda$ in $\RR^n$ is \emph{even} if
$\scal{x}{x} \in 2\ZZ$ for all $x\in\Lambda$;
such a lattice is contained in its dual $\Lambda^*=\{y\in\RR^n\mid\scal{\Lambda}{y}\sseq\ZZ\}$,
and there are integers $\ell\ge1$ such that $\sqrt{\ell}\,\Lambda^*$
is again even ($\ell$ is not unique since, if $\sqrt{\ell}\,\Lambda^*$ is even,
so is $\sqrt{k^2\ell}\,\Lambda^*$ for $k \ge 1$).
For given integers $n$ and $\ell$, denote by
\[
  \Lat_n(\ell)
\]
the class of even lattices $\Lambda$ such that 
$\Lambda ' := \sqrt{\ell}\,\Lambda^*$ is even and
$\det (\Lambda ') = \det (\Lambda)$;
the latter condition implies that $\det (\Lambda) = \ell^{n/2}$.
Clearly, if $\Lambda \in \Lat_n(\ell)$,
then $\Lambda' \in \Lat_n(\ell)$.
A lattice $\Lambda\in\Lat_n(\ell)$ is \emph{$\ell$-modular}
if $\Lambda'$ is equivalent to $\Lambda$.

Now, we describe the modular forms associated to these lattices.
For $\omega$ a nonnegative integer and $\epsilon\in\{+,-\}$, we define
\[
  \varMod^\epsilon_\omega(\ell)
\]
as the space of holomorphic function $f:\HH\to\CC$ 
that verify
\begin{equation*}\tag{$*$}
  \begin{gathered}
    f(z+1) = f(z) \\
    f\Bigl(-\frac{1}{\ell\,z}\Bigr)= \epsilon\, \Bigl(\frac{i}{\sqrt{\ell}\,z}\Bigr)^\omega f(z)
  \end{gathered}
\end{equation*}
for all $z\in\HH$,
and that are holomorphic at infinity, i.e., bounded on $\{z\in\HH\mid \Im z>y_0\}$
for $y_0>0$.
This means that $f$ can be written as
\[
  f(z) = \sum_{m\ge0} a_m q^{2m}, \qquad q=e^{i\pi z}.
\]
Let $f(\infty):= \lim_{\Im z\to\infty} f(z) = a_0$. We define
\[
  \varPod^\epsilon_\omega(\ell) = \{f\in\varMod^\epsilon_\omega(\ell) \mid f(\infty)=0\}.
\]

We give here the classical terminology:
A \emph{modular form} of weight $\omega\geq0$ for a discrete group $\Gamma\subseteq\SL_2(\RR)$
and a character $\chi:\Gamma\to\RR$
is a holomorphic function $f:\HH\to\CC$ that is holomorphic at infinity and that satisfy
\[
  \frac{\chi\bigl(\smatrabcd\bigr)}{(cz+d)^\omega} f\Bigl(\frac{az+b}{cz+d}\Bigr) = f(z) ,
  \qquad \forall \smatrabcd\in\Gamma .
\]
A \emph{parabolic form} is a modular form~$f$ that is zero at infinity,
namely such that $f(\infty)=0$.

For $\ell$ a positive integer, let
$\Gamma_*(\ell)$ be the subgroup of $\SL_2(\RR)$ generated by
\[
  \Gamma_0(\ell) = \Bigl\{ \smatrabcd\in\SL_2(\ZZ) \Bigm| c\equiv 0\bmod \ell \Bigr\}
  \quad \text{and} \quad
  t_\ell = \begin{pmatrix} 0 & 1\big/\sqrt{\ell} \\ -\sqrt{\ell} & 0 \end{pmatrix}.
\]
For any integer $s$,
let $\chi_s$ be the multiplicative character of $\Gamma_*(\ell)$ defined by
\begin{gather*}
  \chi_{s}\Bigl(\smatrabcd\Bigr)=\displaystyle\legendre{(-\ell)^{s}}{d}
                            \quad \text{for $\smatrabcd\in\Gamma_0(\ell)$},\\
  \chi_{s}(t_\ell)=i^s,
\end{gather*}
where $\legendre{\quad}{\quad}$ denotes the Kronecker symbol.
Note that $\chi_s$ depends only of the class of $s$ modulo~$4$.
Denote by $\Mod^+_\omega(\ell)$ the space of modular forms
of weight~$\omega$
for the group $\Gamma_*(\ell)$ and the character $\chi_\omega$, 
and by $\Mod^-_\omega(\ell)$ the space of modular forms
of weight~$\omega$
for the same group and the character $\chi_{\omega+2}$;
denote by $\Pod^{\pm}_\omega(\ell)$ the space of
corresponding parabolic forms.

Note that Equations~$(*)$ say that $\varMod^+_\omega(\ell)$,
respectively $\varMod^-_\omega(\ell)$,
is the space of modular forms of weight~$\omega$ for the subgroup generated by
\[
   T = \begin{pmatrix} 1 & 1 \\ 0 & 1 \end{pmatrix}
   \quad\text{and}\quad
   S_\ell =  \begin{pmatrix} 0 & 1/\sqrt{\ell} \\ -\sqrt{\ell} & 0 \end{pmatrix},
\]
and for the character $\chi_\omega$, respectively $\chi_{\omega+2}$.

\begin{Lemma}\label{IdentificationModularForms}
  For $\ell\in\{1,2,3\}$, the group $\Gamma_*(\ell)$ is generated by
  $T$ and $S_\ell$.
  in particular, for these values of $\ell$, we have
  \[
    \varMod^\epsilon_\omega(\ell) = \Mod^\epsilon_\omega(\ell)
    \quad\text{and}\quad
    \varPod^\epsilon_\omega(\ell) = \Pod^\epsilon_\omega(\ell).
  \]
\end{Lemma}

\begin{proof}
  Let $\ell\in\{1,2,3\}$ and
  $\gamma\in\smatrabcd\in\Gamma_*(\ell)$.
  We show, by induction on $c^2\in\NN$,
  that $\gamma\in\langle T, S_\ell\rangle$,
  where $\langle T, S_\ell\rangle$ denotes the group generated by $T$ and~$S_\ell$.
  
  If $c=0$, then $\gamma=\smatr1&m\\0&1\esmatr=T^m$ for an integer~$m$;
  thus $\gamma\in\langle T, S_\ell\rangle$.
  
  If $c\neq0$, let $k\in\ZZ$ such that $\lvert d+kc\rvert \leq \lvert c/2 \rvert$. We have
  \begin{equation*}
    \gamma\, T^k S_\ell =
      \begin{pmatrix}
        -(b+ka)\sqrt{\ell} & a\big/\sqrt{\ell}\,\null \\
        -(d+kc)\sqrt{\ell} & c\big/\sqrt{\ell}\,\null
      \end{pmatrix},
  \end{equation*}
  and $\bigl(-(d+kc)\sqrt{\ell}\,\bigr)^2 \leq c^2 (\ell/4)  < c^2$.
  By induction, $\gamma\, T^k S_\ell\in\langle T,S_\ell\rangle$, therefore $\gamma\in\langle T,S_\ell\rangle$.
\end{proof}

\begin{Proposition}\label{ThetaSeriesEllLat}
  Let $\Lambda\in\Lat_n(\ell)$,
  where $\ell$ is a positive integer.
  Then
  \begin{align*}
    \Theta_{\Lambda}   + \Theta_{\Lambda'}   &\in \Mod^+_{n/2}(\ell), \\
    \Theta_{\Lambda}   - \Theta_{\Lambda'}   &\in \Pod^-_{n/2}(\ell).
  \end{align*}
  Let moreover $P\in\Har{2h}(\RR^n)$, $h\ge1$. Then,
  if $2h\equiv0\bmod4$,
  \begin{align*}
    \Theta_{\Lambda,P} + \Theta_{\Lambda',P} &\in \Pod^+_{n/2+2h}(\ell), \\
    \Theta_{\Lambda,P} - \Theta_{\Lambda',P} &\in \Pod^-_{n/2+2h}(\ell),
  \end{align*}
  and, if $2h\equiv2\bmod4$,
  \begin{align*}
    \Theta_{\Lambda,P} + \Theta_{\Lambda',P} &\in \Pod^-_{n/2+2h}(\ell), \\
    \Theta_{\Lambda,P} - \Theta_{\Lambda',P} &\in \Pod^+_{n/2+2h}(\ell).
  \end{align*}
\end{Proposition}

\begin{proof}[Partial proof]
  We give the proof only for $\ell\in\{1,2,3\}$. The general case is
  more complicated: see \cite[\S3.1]{Ebeling} together with \cite[Chapter~2]{VenkovMartinet1}.%
  \footnote{Note however that
  our examples of constructions of designs
  in Section~\ref{SectApplications}
  involve lattices in $\Lat_n(\ell)$
  for $\ell\in\{1,2,3\}$ only.}

  It is straightforward that the theta series involved are holomorphic in~$\HH$
  and holomorphic at infinity.
  Now, by Lemma~\ref{IdentificationModularForms}, it is sufficient to
  check that the theta series satisfy Equations~($*$).
  
  As $\Lambda$ and $\Lambda'$ are even, 
  we have clearly $\Theta_{L,P}(z+1)=\Theta_{L,P}(z)$ for $L=\Lambda$~or~$\Lambda'$.
  Now, as a direct consequence of the Poisson Summation Formula---%
  see for example \cite[Prop.~3.1, p.~87]{Ebeling}---we have, for any lattice $L$ of rank~$n$
  and any $P\in\Har{2h}(\RR^n)$,
  \begin{equation*}
    \Theta_{L^*,P}(z) = (\det L)^{1/2} (-1)^h (i/z)^{n/2+2h} \Theta_{L,P}(-1/z).
  \end{equation*}
  From this formula, we deduce, for $L=\Lambda$~or~$\Lambda'$,
  \begin{equation*}
    \Theta_{L',P}(z) = (-1)^h \bigl(i\big/\sqrt{\ell}\,z\bigr)^{n/2+2h} \Theta_{L,P}(-1/\ell z).
  \end{equation*}
  The result follows.
\end{proof} 
  
Let
\begin{gather*}
  \Mod^+(\ell)    := \bigoplus_{\omega\ge0} \Mod^+_\omega(\ell),\\
  \Pod^+(\ell)    := \bigoplus_{\omega>0}   \Pod^+_\omega(\ell),\qquad
  \Pod^-(\ell)    := \bigoplus_{\omega>0}   \Pod^-_\omega(\ell).
\end{gather*}
These are $\Mod^+(\ell)$-modules graded by the weight. For some values of~$\ell$, the structure
of these algebras is known:

\begin{Remark}
  In general, the numerical subscript of modular forms
  given below indicates the \emph{double} of its weight;
  for example $\Delta_{24}\in\Pod_{12}(1)$ is of weight~$12$,
  and $\Theta_{\E_8}\in\Mod_{4}(1)$ is of weight~$4$.
  We have made an exception for the Eisenstein series $E_{(k)}$,
  for which we have retained the traditional notation.
\end{Remark}

\begin{Theorem}\label{ClassifModularForms}
  Let $\ell$ be $1$ or a prime number such that $\ell+1$ divides~$24$.
  Let $k_0$, $k_1$ and $k_2$ given by
  \begin{gather*}
    k_0 =\begin{cases}
            4 & \text{if $\ell=1$,} \\
            2 & \text{if $\ell\equiv 1$ or $2\bmod 4$, and $\ell\neq1$,} \\
            1 & \text{if $\ell\equiv 3\bmod 4$,} \\            
         \end{cases} \\[\jot]
    k_1 = \frac{24}{\ell+1},
    \qquad
    k_2 = k_0+k_1+2.
  \end{gather*}
  Then, we have
  \begin{align*}
    \Mod^+(\ell) &= \CC[\theta_{2k_0},\Delta_{2k_1}], \\
    \Pod^+(\ell) &= \Delta_{2k_1}\, \CC[\theta_{2k_0},\Delta_{2k_1}], \\
    \Pod^-(\ell) &= \Phi_{2k_2}\, \CC[\theta_{2k_0},\Delta_{2k_1}],
  \end{align*}
  where
  \begin{gather*}
    \theta_{2k_0}=\Theta_{L_0}\in\Mod^+_{k_0}(\ell),\text{ with $L_0\in\Lat_{2k_0}(\ell)$,} \\
    \Delta_{2k_1} =\bigl(\eta(z)\,\eta(\ell z)\bigr)^{k_1}\in\Pod^+_{k_1}(\ell), \\
    \Phi_{2k_2}\in\Pod^-_{k_2}(\ell),
  \end{gather*}
  and where $L_0$ is given in Table~\ref{tabSeries1}.
\end{Theorem}

For a proof, see \cite[Section~2]{BachocVenkov} and \cite[Section~3]{Quebbemann}.
Recall that
\[
  \eta(z) = q^{1/12} \prod_{m=1}^{\infty} (1-q^{2m}) = q^{1/12}\bigl(1-q^2-q^4+O(q^{10})\bigr),
                                                                        \qquad q=e^{i\pi z}.
\]
For $\ell=1$, the modular forms of the theorem are $\Delta_{24}=\eta(z)^{24}$,
whose Fourier coefficients are the Ramanujan numbers,
\[
  \Theta_{\E_8} = E_{(4)}, \quad\text{and}\quad \Phi_{36} = \Delta_{24}E_{(6)},
\]
where $E_{(k)}$ is the Eisenstein series of weight~$k$.

Tables \ref{tabSeries1} to \ref{tabSeries3}
give the theta series $\theta_{2k_0}$, $\Delta_{2k_1}$ and $\Phi_{2k_2}$
for all possible values of~$\ell$.
In the third column of Table~\ref{tabSeries1}, the lattice $L_0$ is designated
either by a root system that generates it (e.g., $\E_8$),
or by a quadratic form (e.g., $\smatr2&1\\ 1&4\esmatr$).
In Table~\ref{tabSeries3}, $P_k$ is a suitable element of $\Har{k}(\RR^n)$,
the lattice $L_4$ is the unique (up to isometry) $4$-dimensional lattice of minimum $4$ and determinant $11^2$,
and the lattices $L_0^{(1)}$ and $L_0^{(2)}$ are the two lattices given in the column $L_0$
of Table~\ref{tabSeries1}.
We normalize $\Phi_{2k_2}$ by $\Phi_{2k_2} = q^2 + O(q^4)$.

\begin{table}[!htbp]
  \begin{center}
    \begin{tabular}{|!{\tabstrut} c | c | c | c |}
      \hline
     $\ell$ & $k_0$ & $L_0$ & $\theta_{2k_0}=\Theta_{L_0}$ \\
       \noalign{\hrule}
       $1$  &  $4$  &   $\E_8$                              & $1 + 240q^2 + 2160q^4+ 6720q^6 + O(q^8)$  \\
       $2$  &  $2$  &   $\D_4$                             & $1 +  24q^2 +   24q^4 +  96q^6 + O(q^8)$  \\
       $3$  &  $1$  &   $\A_2\approx\smatr2&1\\ 1&2\esmatr$ & $1 +   6q^2           +   6q^6 + O(q^8)$  \\
       $5$  &  $2$  &   $\A_4$                              & $1 +  20q^2 +   30q^4 +  60q^6 + O(q^8)$  \\
       $7$  &  $1$  &   $\smatr2&1\\ 1&4\esmatr$            & $1 +   2q^2 +    4q^4          + O(q^8)$  \\
      $11$  &  $1$  &   $\smatr2&1\\ 1&6\esmatr$            & $1 +   2q^2           +   4q^6 + O(q^8)$  \\
      $23$  &  $1$  &   $\smatr4&1\\ 1&6\esmatr$ or $\smatr2&1\\1&12\esmatr$
                                                & $1 + 2q^2 + O(q^8)$ or $1 + 2q^4 + 2q^6 + O(q^8)$ \\
      \hline
    \end{tabular}
    \caption{Modular forms $\theta_{2k_0}=\Theta_{L_0}$ of weight~$k_0$}\label{tabSeries1}
  \end{center}
\end{table}

\begin{table}[!htbp]
  \begin{center}
    \begin{tabular}{|!{\tabstrut} c | c | c |}
      \hline
     $\ell$ & $k_1$ & $\Delta_{2k_1} = \bigl(\eta(z)\,\eta(\ell z)\bigr)^{k_1}$ \\
      \hline
       $1$  & $12$  & $q^2 - 24q^4 + 252q^6 - 1472q^8 + O(q^{10})$ \\
       $2$  &  $8$  & $q^2 -  8q^4 +  12q^6 +   64q^8 + O(q^{10})$ \\
       $3$  &  $6$  & $q^2 -  6q^4 +   9q^6 +    4q^8 + O(q^{10})$ \\
       $5$  &  $4$  & $q^2 -  4q^4 +   2q^6 +    8q^8 + O(q^{10})$ \\
       $7$  &  $3$  & $q^2 -  3q^4 +   5q^6 -    7q^8 + O(q^{10})$ \\
      $11$  &  $2$  & $q^2 -  2q^4 -    q^6 +    2q^8 + O(q^{10})$ \\
      $23$  &  $1$  & $q^2 -   q^4 -    q^6           + O(q^{10})$ \\
      \hline
    \end{tabular}
    \caption{Modular forms $\Delta_{2k_1} = \bigl(\eta(z)\,\eta(\ell z)\bigr)^{k_1}$~of weight $k_1$}\label{tabSeries2}
  \end{center}
\end{table}
\begin{table}[!htbp]
  \begin{center}
    \begin{tabular}{|!{\vrule depth 7pt height 13pt width 0pt} c | c | c | c |}
      \hline
     $\ell$ & $k_2$ & $\Phi_{2k_2}$ (definition) & $\Phi_{2k_2}$ (expansion) \\
      \hline
       $1$  & $18$  & $\Theta_{\E_8,P_{14}}$                     & \small $q^2 - 528q^4 -4284q^6 +147712q^8 + O(q^{10})$ \\
       $2$  & $12$  & $\Theta_{\D_4^2,P_8} - \Theta_{(\D_4^2)',P_8}$ & $q^2 -  88q^4 + 252q^6 +    64q^8 + O(q^{10})$ \\
       $3$  &  $9$  & $\Theta_{\A_2^3,P_6} + \Theta_{(\A_2^3)',P_6}$ &  $q^2 -  14q^4 +  48q^6 +    68q^8 + O(q^{10})$ \\
       $5$  &  $8$  & $\Theta_{L_0,P_6}    + \Theta_{L_0',P_6}$      & $q^2 -  14q^4 -  48q^6 +    68q^8 + O(q^{10})$ \\
       $7$  &  $7$  & $\Theta_{L_0^2,P_4}  - \Theta_{(L_0^2)',P_4}$  & $q^2 -  10q^4 -  14q^6 +    68q^8 + O(q^{10})$ \\
      $11$  &  $5$  & $\Theta_{L_0\perp L_4,P_2} + \Theta_{(L_0\perp L_4)',P_2}$
                                                                     & $q^2 -   6q^4 -   3q^6 -    14q^8 + O(q^{10})$ \\
      $23$  &  $4$  & \small$\Theta_{L_0^{(1)}\perp L_0^{(2)},P_2} + \Theta_{(L_0^{(1)}\perp L_0^{(2)})',P_2}$
                                                                     & $q^2 -   2q^4 -   5q^6 -     4q^8 + O(q^{10})$ \\
      \hline
    \end{tabular}
    \caption{Modular forms $\Phi_{2k_2}$ of weight~$k_2$}\label{tabSeries3}
  \end{center}
\end{table}

\section{Construction of cubature formulas using shells of lattices}\label{SectConstruction}

Let $G$ be a finite subgroup of $\Orth(n)$,
and let $X_1, X_2, \dots, X_r$ be nonempty finite subsets of $\Se^{n-1}$
such that $G\,X_j=X_j$.
Consider the set
\[
  X := \bigcup_{j=1}^{r} X_j \quad\sseq\quad \Se^{n-1}.
\]
The aim is to find numbers $W_1, W_2, \dots, W_r\in\RR_{>0}$
such that $(X,W)$ is a cubature formula of high strength,
where $W$ is the weight function given by
\[
  W(x) = \sum_{j\,:\, x\in X_j} W_j.
\]
By Proposition~\ref{CharactHarmonic},
$(X,W)$ is of strength~$t$ if and only if
the following conditions on
$W_1,\dots,W_r$ are satisfied:
\begin{gather*}
  \sum_{j=1}^{r} W_j \card{X_j} = 1,  \\
  \sum_{j=1}^{r} W_j\, P(X_j) = 0 ,\qquad \forall P\in\Har{k},\ 1\le k \le t,
\end{gather*}
where
\[
  P(X_j):=\sum_{x\in X_j} P(x).
\]
The following statement, together with Proposition~\ref{CharactHarmonicInv},
is fundamental in our analysis.

\begin{Lemma}\label{LinearEquations}
  Let $G$ be a finite subgroup of $\Orth(n)$, and
  let $E\sseq\RR^n$ be a discrete set that is invariant by~$G$.
  Let $m_1, m_2, \dots, m_r$ be pairwise distinct positive numbers, let
  \[
    X_j = \frac{1}{\sqrt{m_j}} E_{m_j}\quad \sseq\quad \Se^{n-1},
    \qquad j=1,\dots,r,
  \]
  where $E_{m} = \{ x\in E \mid \scal{x}{x}=m\}$.
  Let $X = \bigcup_{j=1}^{r}X_j$ and $W_1, \dots, W_r>0$, and
  let $k$ be a positive integer.
  
  Assume that there exist formal series
  \[
    \Theta_i = \sum_{m>0} a_i(m)\,q^m,\qquad i=1,\dots,N
  \]
  such that, for every $P\in\Har{k}(\RR^n)^G$, the theta series $\Theta_{E,P}$ is of the form
  \[
    \Theta_{E,P} = c_1(P)\,\Theta_1 + c_2(P)\,\Theta_2 + \dots +c_N(P)\,\Theta_N
  \]
  for some $c_i(P)\in\CC$.
  
  Then, for the condition
  \[
    \sum_{j=1}^{r} W_jP(X_j) = 0 \quad \forall P\in\Har{k}(\RR^n)^G
  \]
  to hold, it suffices that
  \[\tag{$\sharp$}
    \sum_{j=1}^r \frac{a_i(m_j)}{m_j^{k/2}} W_j = 0, \quad i=1,\dots,N.
  \]
\end{Lemma}

\begin{proof}
Assume that Equation~$(\sharp)$ holds,
and let $P\in\Har{k}(\RR^n)^G$. By hypothesis, there exist
$c_1,\cdots,c_N\in\CC$ such that
\[
  \sum_{m>0} P(E_m) q^m = \Theta_{E,P} = \sum_{i=1}^{N}c_i\,\Theta_i =  \sum_{m>0} \biggl( \sum_{i=1}^{N} c_i\, a_i(m) \biggr)  q^m.
\]
Since $P$ is homogeneous of degree~$k$, we have $P\bigl( (1/\!\sqrt{m})\,E_m\bigr) = m^{-k/2} P(E_m)$.
Therefore,
\[
  P(X_j)  = m_j^{-k/2} P(E_{m_j}) = \sum_{i=1}^{N} c_i \frac{a_i(m_j)}{m_j^{k/2}},
\]
hence
\[
  \sum_{j=1}^{r} W_j P(X_j) = \sum_{i=1}^{N} c_i \sum_{j=1}^{r} \frac{a_i(m_j)}{m_j^{k/2}}W_j = 0.
  \rlap{\quad\null\qedhere}
\]
\end{proof}

\section{Applications}\label{SectApplications}

If not otherwise specified, Lemma~\ref{LinearEquations} is applied with
the trivial group $G=\{\mathrm{id}\}$.
In this section, we do not give the details for all lattices we have considered,
but we have made a selection which reflects most situations that may occur.
All cubature formulas we have obtained are listed
in Section~\ref{SectResults} below.

\subsection{The Leech lattice}

Let $\Lambda$ be the Leech lattice, namely the
unique (up to isometry)
even unimodular lattice of dimension~$24$
and of minimum~$4$.
Let $P\in\Har{2h}(\RR^n)$.
From Proposition~\ref{ThetaSeriesEllLat} and Theorem~\ref{ClassifModularForms},
we deduce (using $\Lambda'=\Lambda$):
\[
  \Theta_{\Lambda,P} =
    \begin{cases}
      \Theta_{\E_8}^3 - 720\Delta_{24}     & \text{if $P=1$,} \\
      0                                    & \text{if $P\in\Har{2h}(\RR^{24})$, $2h=2,4,6,8,10,14$,} \\
      c_1(P)\,\Delta_{24}^2                & \text{if $P\in\Har{12}(\RR^{24})$,} \\
      c_2(P)\,\Theta_{\E_8}\,\Delta_{24}^2 & \text{if $P\in\Har{16}(\RR^{24})$,} \\
   \end{cases}
\]
where $c_1$ and $c_2$ are linear forms on $\Har{12}(\RR^{24})$ and $\Har{16}(\RR^{24})$
respectively.
The first coefficients of these theta series are
\begin{gather*}
      \Theta_{\E_8}^3 - 720\,\Delta_{24} = 1 + 196\,560\,q^4 + 16\,752\,960\,q^6 + O(q^8), \\
      \Delta_{24}^2 = q^4 -  48q^6 + O(q^8), \\
      \Theta_{\E_8}\,\Delta_{24}^2 = q^4 + 192q^6 + O(q^8).
\end{gather*}
Since $\Theta_{\Lambda,P}=0$ for $P\in\Har{2h}(\RR^{24})$, $2\leq 2h\leq 10$,
it follows from Proposition~\ref{CharactHarmonic} that every shell of the Leech lattice is a spherical $11$-design.

Now, we want to find cubature formulas of higher strengths
using the construction of Section~\ref{SectConstruction}.
Consider the set
\[
  X = X_1 \cup X_2, \quad \text{where}\quad X_1 = \frac{1}{\sqrt{4}}\Lambda_4,\ X_2 = \frac{1}{\sqrt{6}}\Lambda_6.
\]
The two sets $X_1$ and $X_2$ are disjoint.
Indeed, if we had $x\in X_1\cap X_2$,
we would have $\sqrt{4}\,x\in \Lambda_4$ and $\sqrt{6}\,x\in\Lambda_6$;
but $\scal{\sqrt{4}\,x}{\sqrt{6}\,x} = 2\sqrt{6}\scal{x}{x} = 2\sqrt{6}$,
which is impossible since $\Lambda$ is integral.

We want to find out the numbers $W_1$ and $W_2$
such that $(X,W)$
is a cubature of strength~$15$
(same notation as in Section~\ref{SectConstruction}).
For this, we have to fulfill the conditions
\begin{gather*}
  \sum_{x\in X} W(x)=1, \\
  \sum_{x\in X} W(x)\,P(x)=0 \text{ for } P=\Har{12}(\RR^n).
\end{gather*}
The first condition is equivalent to $\card{X_1}W_1 + \card{X_2}W_2=1$, that is
\[
  196\,560\, W_1 +  16\,752\,960\, W_2 = 1.
\]
By Lemma~\ref{LinearEquations} the second condition is equivalent to
\[
    \frac{1}{4^{6}}\,W_1 - \frac{48}{6^{6}}\,W_2 = 0,
\]
or $W_2 = (3/4)^5 W_1$.
The solution of this set of two linear equations is
\[
  W_1 \approx 2.394\times 10^{-7},\qquad W_2\approx 0.568\times 10^{-7}.
\]
So we get a cubature formula of size $196\,560 + 16\,752\,960 = 16\,949\,520$.

We can obtain in the same way a cubature formula of strength~$17$ by using the shells
of norm~$4$,~$6$, and~$8$, and a cubature formula of strength~$19$ by using the shells
of norm~$4$,~$6$, $8$, and~$10$ of the Leech lattice.

\subsection{The Korkine-Zolotareff lattice}

Let $\Lambda$ be the Korkine-Zolotareff lattice, namely the unique (up to isometry) even unimodular lattice of dimension~$8$.
It is the lattice generated by the root system of type~$\E_8$.
From Proposition~\ref{ThetaSeriesEllLat} and Theorem~\ref{ClassifModularForms},
we deduce (using $\Lambda'=\Lambda$):
\[
  \Theta_{\Lambda,P} =
    \begin{cases}
      \Theta_{\E_8}                        & \text{if $P=1$,} \\
      0                                    & \text{if $P\in\Har{2h}(\RR^{8})$, $2h=2,4,6,10$,} \\
      c_1(P)\,\Delta_{24}                  & \text{if $P\in\Har{8}(\RR^{8})$,} \\
      c_3(P)\,\Theta_{\E_8} \Delta_{24}    & \text{if $P\in\Har{12}(\RR^{8})$,} \\
      c_2(P)\,\Phi_{36}                    & \text{if $P\in\Har{14}(\RR^{8})$.} 
   \end{cases}
\]
So, by Proposition~\ref{CharactHarmonic}, any shell of the Korkine-Zolotareff lattice is a spherical $7$-design,
and we can obtain a cubature formula of stength~$11$ by combining the shells of norms $2$~and~$4$.

But if we try to obtain a cubature formula of strength~$13$ of nodes
\[
  X = X_1 \cup X_2\cup X_3, \qquad\text{where}\quad
  X_1 = \frac{1}{\sqrt{2}}\Lambda_2,\ X_2 = \frac{1}{\sqrt{4}}\Lambda_4,\ X_3 = \frac{1}{\sqrt{6}}\Lambda_6,
\]
we obtain the weights
\[
  W_1 \approx -0.744\times 10^{-4}, \quad W_2 \approx 1.587\times 10^{-4}, \quad  W_3 \approx 1.005\times 10^{-4},
\]
with $W_1$ negative.
A similar problem occurs if we try to use the shells of norm $2$,~$4$ and~$8$.
Therefore, in order to get a true cubature formula with positive weights, we have to use
\[
  X = X_1 \cup X_2\cup X_3, \qquad X_1 = \frac{1}{\sqrt{2}}\Lambda_2,\ X_2 = \frac{1}{\sqrt{6}}\Lambda_6,\ X_3 = \frac{1}{\sqrt{8}}\Lambda_8.
\]
We obtain the weights
\[
  W_1 \approx 0.792\times 10^{-4}, \quad W_2 \approx 0.457\times 10^{-4}, \quad  W_3 \approx 0.385\times 10^{-4}.
\]
Here, we have $X_1\sseq X_3$, since $2x\in\Lambda_8$ if $x\in\Lambda_2$.
Therefore, we obtain a cubature formula
of strength~$13$ of size
\[
  \card{X} = \card{X_2}+\card{X_3} = 6\,720 + 17\,520 = 24\,240
\]
and of weights
\[
  W(x) =
  \begin{cases}
    W_2     \approx 0.457\times 10^{-4}  & \text{if $x\in X_2$,} \\
    W_3     \approx 0.385\times 10^{-4} & \text{if $x\in X_3\setminus X_1$,} \\
    W_1+W_3 \approx 1.177\times 10^{-4} & \text{if $x\in X_1$.}    
  \end{cases}
\]

\subsection{The Barnes-Wall lattice of dimension~16}
Let $\Lambda$ be the Barnes-Wall lattice of dimension~$16$. It is a $2$-modular lattice, that is
$\Lambda':=\sqrt{2}\Lambda$ is equivalent to~$\Lambda$, and it is of minimum~$4$.
From Proposition~\ref{ThetaSeriesEllLat} and Theorem~\ref{ClassifModularForms},
we have:
\begin{gather*}
  \\[\jot]
  \Theta_{\Lambda,P}+\Theta_{\Lambda',P} =
    \begin{cases}
      2\,\Theta_{\D_4}^4-192\,\Delta_{16}           & \text{if $P=1$,} \\
      0                                             & \text{if $P\in\Har{2h}(\RR^{16})$, $2h=2,4,6,10$,} \\
      c_1(P)\,\Delta_{16}^2                         & \text{if $P\in\Har{8}(\RR^{16})$,} \\
      c_2(P)\,\Delta_{16}^2\,\Theta_{\D_4}^2        & \text{if $P\in\Har{12}(\RR^{16})$,} \\
      c_3(P)\,\Phi_{24}\,\Delta_{16}\,\Theta_{\D_4} & \text{if $P\in\Har{14}(\RR^{16})$,} \\
    \end{cases}
  \\[\jot]
  \Theta_{\Lambda,P} - \Theta_{\Lambda',P} =
    \begin{cases}
      0                                             & \text{if $P\in\Har{2h}(\RR^{16})$, $2h=0,2,4,6,8,12$,} \\
      c_4(P)\,\Delta_{16}^2\,\Theta_{\D_4}          & \text{if $P\in\Har{10}(\RR^{16})$,} \\
      c_5(P)\,\Phi_{24}\,\Delta_{16}                & \text{if $P\in\Har{12}(\RR^{16})$,} \\
    \end{cases}
\end{gather*}  
and therefore
\[
  \Theta_{\Lambda,P} =
    \begin{cases}
      \Theta_{\D_4}^4-96\,\Delta_{16}       & \text{if $P=1$,} \\
      0                                    & \text{if $P\in\Har{2h}(\RR^{16})$, $2h=2,4,6$,} \\
      c_6(P)\,\Delta_{16}^2                & \text{if $P\in\Har{8}(\RR^{16})$,} \\
      c_7(P)\,\Delta_{16}^2\Theta_{\D_4}   & \text{if $P\in\Har{10}(\RR^{16})$.} \\
    \end{cases}
\]
Every shell of $\Lambda$ is a spherical $7$-design, and we obtain
a cubature formula of strength~$9$ by combining two different shells of~$\Lambda$,
for example $\Lambda_4$ and $\Lambda_6$.
We also obtain a cubature formulas of strength~$11$
of support $\frac{1}{\sqrt4}(\Lambda\cup\nobreak\Lambda')_4 \cup \frac{1}{\sqrt6}(\Lambda\cup\nobreak\Lambda')_6$,
and a cubature formulas of strength~$13$
of support $\frac{1}{\sqrt4}(\Lambda\cup\nobreak\Lambda')_4 \cup \frac{1}{\sqrt6}(\Lambda\cup\nobreak\Lambda')_6
\cup \frac{1}{\sqrt{10}}(\Lambda\cup\nobreak\Lambda')_{10}$.

\subsection{The $\D_4$ root lattice}

Let $\Lambda$ be the lattice generated by the root system $\D_4$.
It is a $2$-modular lattice. Here, we will apply Lemma~\ref{LinearEquations}
with 
\[
  G=\Aut(\Lambda\cup\Lambda').
\]
We now describe $\Lambda$, $\Lambda'$ and $G$. We have
\begin{gather*}
  \Lambda  = \bigl\{ x=(x_1,x_2,x_3,x_4) \bigm| x_i\in \ZZ,\ x_1+x_2+x_3+x_4\in 2\ZZ \bigr\}, \\
  \Lambda' = \bigl\{ x=(x_1,x_2,x_3,x_4) \bigm| x_i\in{\textstyle\frac{1}{\sqrt{2}}}\ZZ,\ x_i-x_j\in\sqrt{2}\ZZ \bigr\} .
\end{gather*}
The group $\Aut(\Lambda)$
contains (actually, is) a reflection group,
denoted by $\Weyl(\F_4)$,
which is generated by the 24~reflections $x\mapsto x - 2\frac{\scal x\alpha}{\scal\alpha\alpha}\alpha$,
where $\alpha$ is one of the following vectors:
\[
  (\pm1,0,0,0), \quad (\pm1,\pm1,0,0), \quad (\pm1,\pm1,\pm1,\pm1)
\]
(all choices of signs and all permutations of coordinates).
The invariant polynomials of this group are known \cite[\S7.4, pp.~217--218]{Smith};
they form a polynomial $\CC$\nobreakdash-algebra of basis $\{h_2,h_6,h_8,h_{12}\}$, where
\begin{gather*}
  h_2 = \sym(x_1^2), \qquad
  h_6 = \sym(x_1^4 x_2^2) - 3 \sym(x_1^2 x_2^2 x_3^2), \\[\jot]
  h_8 = \sym(x_1^8) + 14\sym(x_1^4 x_2^4) + 168\,x_1^2 x_2 ^2 x_ 3^2 x_4^2, \\[\jot]
  \begin{split}
     h_{12} = \sym(x_1^{12}) + 22\sym(x_1^6 x_2^6) +165\sym(x_1 ^4 x_2^4 x_3^4) \qquad \\
          \qquad {}+ 330 \sym(x_1^6 x_2 ^2 x_3 ^2 x_4^2) + 330\sym(x_1^4 x_2^4 x_3^2 x_4^2),
  \end{split}
\end{gather*}
and where
\[
  \sym(p) =  \frac{1}{\card{(\Sym_4)_p}}\sum_{\sigma\in{\Sym_4}} \sigma\cdot p,
  \qquad (\Sym_4)_p:=\{\sigma\in\Sym_4\mid \sigma\cdot p=p\},
\]
where the symmetric group $\Sym_4$ acts by permuting the coordinates.
For our purpose, it is better to choose the basis
\begin{gather*}
   H_2:=h_2, \qquad
   H_6:=8h_6-h_2^3, \qquad
   H_8:=10h_8-7h_2^4, \\
   H_{12} := 64 h_{12} - 55 h_8 h_2^2 - 176 h_6^2 + 220 h_6 h_2^3 - 11 h_2^6.
\end{gather*}
The orthonormal transformation
\[
  T:=\begin{pmatrix}
       -\sqrt{2}/2 &  \sqrt{2}/2 &           0 &           0 \\
       -\sqrt{2}/2 & -\sqrt{2}/2 &           0 &           0 \\
                 0 &           0 &  \sqrt{2}/2 &  \sqrt{2}/2  \\
                 0 &           0 & -\sqrt{2}/2 &  \sqrt{2}/2
      \end{pmatrix}
\]
exchanges $\Lambda$ and $\Lambda'$. So, $G=\Aut(\Lambda\cup\Lambda')$ is generated by
$\Aut(\Lambda)$ and $T$.
The action of $T$ on the polynomials $H_k$ is:
\[
  TH_2 = H_2,\qquad TH_6=-H_6,\qquad TH_8=H_8, \qquad TH_{12} = -H_{12}.
\]
Therefore, the polynomials that are invariant by~$G$ are linearly generated by the polynomials
\[
  H_2^\alpha H_6^\beta H_8^\gamma H_{12}^\delta, \qquad
  \beta+\delta\equiv0\bmod2.
\]
Let $\Pol{k}(\RR^n)$ be the space of homogeneous polynomial functions on~$\RR^n$ of degree~$k$.
Since $\Har{k}(\RR^n)$ is the kernel of the surjective $G$-equivariant Laplacian
$\Delta:\Pol{k}(\RR^n)\to\Pol{k-2}(\RR^n)$, we have
\[
  d_k^G := \dim \Har{k}(\RR^n)^G = \dim \Pol{k}(\RR^n)^G - \dim \Pol{k-2}(\RR^n)^G.
\]
We can compute
\[
  \sum_{k\ge0} d_k^G X^k = \frac{1+X^{18}}{(1-X^8)(1-X^{12})(1-X^{24})} = 1 + X^8 + X^{12} + X^{16} + O(X^{18}) .
\]
So, we have by Proposition~\ref{ThetaSeriesEllLat},
\[
  \Theta_{\Lambda,P}+\Theta_{\Lambda',P} =
    \begin{cases}
      2\,\Theta_{\D_4}                       & \text{if $P=1$,} \\
      0                                      & \text{if $P\in\Har{2h}(\RR^{8})^G=\{0\}$, $2h=2,4,6,10,14$,} \\
      c_1(P)\,\Delta_{16}\,\Theta_{\D_4}     & \text{if $P\in\Har{8}(\RR^{8})^G$,} \\
      c_2(P)\,\Delta_{16}^2\,\Theta_{\D_4}^3 & \text{if $P\in\Har{12}(\RR^{8})^G$.} \\
    \end{cases}
\]
Thus, using Lemma~\ref{LinearEquations}
together with Proposition~\ref{CharactHarmonicInv},
we obtain cubature formulas of strengths $11$ and~$15$
by combining two or three shells of $\Lambda\cup \Lambda'$.
Note however that certain combinations of shells,
such as $\frac{1}{\sqrt2}(\Lambda\cup \Lambda')_2\cup \frac{1}{\sqrt4}(\Lambda\cup \Lambda')_4$
are not possible,
because Lemma~\ref{LinearEquations} can provide a degenarate linear equations system.
(In the present example, we have indeed
$\frac{1}{\sqrt2}(\Lambda\cup \Lambda')_2 = \frac{1}{\sqrt4}(\Lambda\cup \Lambda')_4$.)

\subsection{The shorter Leech lattice}

Let $\Lambda$ be the unique (up to isometry) odd unimodular integral lattice of dimension~23 and of minimum~3;
it is the \emph{shorter Leech lattice}.
As for even ones,
the theta series of odd unimodular lattices can be computed: see \cite{hpv}.
So, we can apply our method to construct cubature formulas from the shells of~$\Lambda$.

\section{Cubature formulas obtained by our construction}\label{SectResults}

The cubature formulas we obtain by our method
are summarized in Tables \ref{tabCF4}~to~\ref{tabCF24}.
These tables involved the lattices listed in Tables~\ref{tabLatticesEven}~and~\ref{tabLatticesOdd}.

In Tables \ref{tabCF4} to \ref{tabCF24}, the column ``bound'' indicates the lower bound
of the size
of cubature formulas of strength~$t$ on $\Se^{n-1}$
obtained by Linear Programming.
The indication in parentheses has the following meaning:
\begin{myitemize}
\item (LP$d$) refers to the estimate of the Linear Programming bound
  obtained by the procedure described in Section~\ref{SectBounds},
  paragraph ``Numerical estimate of the Linear Programming bound'',  
  where $2d$ is the degree of the polynomial used for the estimate.
  The calculations have been performed with ``Maple'' software,
  using the ``Optimization'' package.
  Because of limitations in time and memory of our computer,
  we have not always chosen an optimal~$d$,
  in which case we write (LP${\ge}d$);
  similarly, when the bound we have found
  can almost certainly be improved, we put the sign $\ge$ before it.
\item (T) refers to the Delsarte or (Fisher-type) bound;
  it is the same as (LP$d$) when $2d+1=t$ (the strength).
\item When it is known that tight spherical designs
  with the corresponding parameters $(n,t)$ do not exist,
  we note (D) to indicate the Delsarte bound increased by~1
\end{myitemize}

The columns ``set'' and ``shells'' indicate respectively the set $E$ and the shells $m_j$
used in the construction of Section~\ref{SectConstruction},
according to the notation of Lemma~\ref{LinearEquations}.
The column ``size'' gives the size of the cubature formula so constructed.

In the case $n=4$, we have also listed in Table~\ref{tabCF4b}
the cubature formulas of smallest known sizes for fixed values of~$t$;

Most spherical designs with $n\le24$ described in \cite[Table~3, p.~108]{BachocVenkov}
appear also in our tables.

\begin{table}[H]
  \begin{center}
    \begin{tabular}{|!{\tabstrut} c | c | c | c |}
      \hline
       $n$ & $\ell$ & notation    & name \\\hline\hline
       $4$ & $2$ & $D_4$          & $\D_4$ root lattice \\\hline
       $8$ & $1$ & $\KZ_8$        & Korkine-Zolotareff \\\hline
      $12$ & $3$ & $K_{12}$       & Coxeter-Todd \\\hline
      $14$ & $3$ & $Q_{14}$       & Quebbemann \\\hline
      $16$ & $2$ & $\BW_{16}$     & Barnes-Wall \\\hline
      $20$ & $2$ & $N_{20}^{(1)}, N_{20}^{(2)}, N_{20}^{(3)}$
                       & Nebe \\\hline
      $24$ & $3$ &$N_{24}$        & Nebe \\\hline
      $24$ & $1$ & $\Lambda_{24}$ & Leech \\\hline
    \end{tabular}
    \caption{Even $\ell$-modular lattices involved in Tables~\ref{tabCF4}~to~\ref{tabCF24}}\label{tabLatticesEven}
  \end{center}
\end{table}

\begin{table}[H]
  \begin{center}
    \begin{tabular}{|!{\tabstrut} c | c | c |}
      \hline
       $n$ & notation & name \\\hline\hline
      $23$ & $O_{23}$ & shorter Leech lattice \\\hline
    \end{tabular}
    \caption{Odd unimodular lattice involved in Table~\ref{tabCF23}}\label{tabLatticesOdd}
  \end{center}
\end{table}

\begin{table}[H]
 \begin{center}
   \begin{tabular}{|!{\tabstrut} c | c | c | c | c |}
     \hline
     strength & set            & shells    & size  & bound        \\\hline\hline
     $5$      & $D_4$          & $2$       &  $24$ &  $21$ (D)    \\\hline
     $7$      & $D_4\cup D_4'$ & $2$       &  $48$ &  $42$ (LP5)  \\\hline  
     $11$     & $D_4\cup D_4'$ & $2,6$     & $240$ & $120$ (LP11) \\\hline 
     $15$     & $D_4\cup D_4'$ & $2,6, 10$ & $528$ & $267$ (LP22) \\\hline 
   \end{tabular}
  \caption{Cubature formulas for $n=4$ ---see also Table\ref{tabCF4b}}\label{tabCF4}
 \end{center}
\end{table}
 
\begin{table}[H]
 \begin{center}
   \begin{tabular}{|!{\tabstrut} c | l | c | c |}
     \hline
     strength &\null\hfill best known cubature formula\hfill\null & size & bound \\\hline\hline
     $3$  & root system $\A_1^4$             &  $8$ &      $8$ (T)         \\\hline
     $5$  & root system $\D_4$               & $24$ &     $21$ (D)         \\\hline
     $7$  & \cite{HardinSloane}              & $46$ &     $42$ (LP5)       \\\hline 
     $9$  & announced in \cite{HardinSloane} & $86$ &     $74$ (LP7)       \\\hline 
     $11$ & vertices of the 600-cell         &$120$ &    $120$ (LP11)      \\\hline 
     $19$ & $\vcenter{\kern2pt\hbox{\strut vertices of the 600-cell and of}
                      \hbox{its dual, the $120$-cell \cite{Salihov}}\kern2pt}$
                                             &$720$ & $500$ (LP$\ge$18) \\\hline 
  \end{tabular}
  \caption{Best known cubature formulas for $n=4$}\label{tabCF4b}
 \end{center}
\end{table}

\begin{table}[H]
 \begin{center}
   \begin{tabular}{|!{\tabstrut} c | c | c | c | c |}
     \hline
     strength & set & shells & size & bound \\\hline\hline
     $7$      & $\KZ_{8}$ & $2$        &     $240$ &       $240$ (T)     \\\hline
     $11$     & $\KZ_{8}$ & $2,4$      &  $2\,400$ &    $1\,856$ (LP8)   \\\hline 
     $13$     & $\KZ_{8}$ & $2,6,8$    & $24\,240$ &    $4\,361$ (LP12)  \\\hline 
     $15$     & $\KZ_{8}$ & $2,4,6,8$  & $26\,400$ &    $9\,190$ (LP16)  \\\hline 
   \end{tabular}  
  \caption{Cubature formulas for $n=8$}\label{tabCF8}
 \end{center}
\end{table}

\begin{table}[H]
 \begin{center}
   \begin{tabular}{|!{\tabstrut} c | c | c | c | c |}
     \hline
     strength & set & shells & size & bound \\\hline\hline
     $5$   & $K_{12}$             & $4$      &       $756$ &     $157$ (D)   \\\hline 
     $7$   & $K_{12}\cup K_{12}'$ & $4$      &    $1\,512$ &     $729$ (D)   \\\hline 
     $9$   & $K_{12}\cup K_{12}'$ & $6$      &    $8\,064$ &  $2\,940$ (LP6) \\\hline 
     $11$  & $K_{12}\cup K_{12}'$ & $4,6,8$  &   $50\,400$ & $10\,604$ (LP7) \\\hline 
   \end{tabular}  
  \caption{Cubature formulas for $n=12$}\label{tabCF12}
 \end{center}
\end{table}

\begin{table}[H]
 \begin{center}
   \begin{tabular}{|!{\tabstrut} c | c | c | c | c |}
     \hline
     strength & set & shells & size & bound \\\hline\hline
     $5$   & $Q_{14}$             & $4$      &      $756$ &       $211$ (D)   \\\hline
     $7$   & $Q_{14}\cup Q_{14}'$ & $4$      &   $1\,512$ &    $1\,121$ (D)   \\\hline
     $9$   & $Q_{14}\cup Q_{14}'$ & $4,8$    &  $89\,964$ &    $4\,902$ (LP6) \\\hline 
     $11$  & $Q_{14}\cup Q_{14}'$ & $4,6,8$  & $107\,436$ &   $20\,817$ (LP7) \\\hline
   \end{tabular}  
  \caption{Cubature formulas for $n=14$}\label{tabCF14}
 \end{center}
\end{table}

\begin{table}[H]
 \begin{center}
   \begin{tabular}{|!{\tabstrut} c | c | c | c | c |}
     \hline
     strength & set & shells & size & bound \\\hline\hline
     $7$   & $\BW_{16}$              & $4$       &     $4\,320$ &     $1\,633$ (D)   \\\hline
     $9$   & $\BW_{16}$              & $4,6$     &    $65\,560$ &     $7\,753$ (D)   \\\hline
     $11$  & $\BW_{16}\cup\BW_{16}'$ & $4,6$     &   $131\,520$ &    $37\,166$ (LP7) \\\hline
     $13$  & $\BW_{16}\cup\BW_{16}'$ & $4,6,10$  &$4\,555\,200$ &$\ge146\,153$ (LP8) \\\hline
   \end{tabular}  
  \caption{Cubature formulas for $n=16$}\label{tabCF16}
 \end{center}
\end{table}

\begin{table}[H]
 \begin{center}
   \begin{tabular}{|!{\tabstrut} c | c | c | c | c |}
     \hline
     strength & set & shells & size & bound \\\hline\hline
     $5$   & $N_{20}$             & $4$       &     $3\,960$ &        $421$ (D)   \\\hline
     $7$   & $N_{20}\cup N_{20}'$ & $4$       &     $7\,920$ &     $3\,081$ (D)   \\\hline
     $9$   & $N_{20}\cup N_{20}'$ & $4,6$     &   $345\,840$ &    $17\,711$ (D)   \\\hline
     $11$  & $N_{20}\cup N_{20}'$ & $4,6,8$   &$4\,527\,600$ &   $95\,309$ (LP7)  \\\hline 
   \end{tabular}  
   \endgraf\smallskip
   In this table, $N_{20}$ denotes any of the three lattices $N_{20}^{(1)}$, $N_{20}^{(2)}$, and~$N_{20}^{(3)}$.
  \caption{Cubature formulas for $n=20$}\label{tabCF20}
 \end{center}
\end{table}

\begin{table}[H]
 \begin{center}
   \begin{tabular}{|!{\tabstrut} c | c | c | c | c |}
     \hline
     strength & set & shells & size & bound \\\hline\hline
     $7$      & $O_{23}$  & $3$     &     $4\,600$ &     $4\,600$ (T)   \\\hline
     $9$      & $O_{23}$  & $3,5$   &   $958\,458$ &    $29\,901$ (D)   \\\hline
     $11$     & $O_{23}$  & $3,4,6$ &$6\,574\,550$ &   $166\,808$ (LP7) \\\hline  
   \end{tabular}  
  \caption{Cubature formulas for $n=23$}\label{tabCF23}
 \end{center}
\end{table}

\begin{table}[H]
 \begin{center}
   \begin{tabular}{|!{\tabstrut} c | c | c | c | c |}
     \hline
     strength & set & shells & size & bound \\\hline\hline
     $5$  & $N_{24}$             & $6$        &         $26\,208$ &              $601$ (D)     \\\hline
     $7$  & $N_{24}\cup N_{24}'$ & $6$        &         $52\,416$ &           $5\,201$ (D)     \\\hline
     $11$ & $\Lambda_{24}$       & $4$        &        $196\,560$ &         $196\,560$ (T)     \\\hline
     $15$ & $\Lambda_{24}$       & $4,6$      &    $16\,969\,680$ & $\ge  6\,179\,991$ (LP9)   \\\hline 
     $17$ & $\Lambda_{24}$       & $4,6,8$    &   $415\,003\,680$ & $\ge 27\,131\,261$ (LP13)  \\\hline 
     $19$ & $\Lambda_{24}$       & $4,6,8,10$ &$5\,044\,384\,800$ & $\ge116\,303\,274$ (LP15)  \\\hline 
   \end{tabular}  
  \caption{Cubature formulas for $n=24$}\label{tabCF24}
 \end{center}
\end{table}


\begin{thebibliography}{HarPac04b}
\bibitem[NikNik02]{NikovaNikov}
    S.~Nikova, V.~Nikov,
    \emph{Improvement of the Delsarte bound for $\tau$-designs when it is not the best bound possible},
     Des.\ Codes Cryptogr.\ \textbf{28} (2) (2003) 201--222.
\end{thebibliography}
\end{document}